\documentclass{article}
\usepackage{amsmath,amsxtra,amssymb,latexsym, amscd}

\usepackage{graphicx}
\begin{document}
\newtheorem{thm}{Theorem}
\newtheorem{cor}[thm]{Corollary}
\newtheorem{lem}[thm]{Lemma}
\newtheorem{prop}[section]{Proposition}
\newtheorem{defn}[thm]{Definition}
\newtheorem{rem}[thm]{Remark}
\numberwithin{equation}{section}
\newcommand{\norm}[1]{\left\Vert#1\right\Vert}
\newcommand{\abs}[1]{\left\vert#1\right\vert}
\newcommand{\set}[1]{\left\{#1\right\}}
\newcommand{\Real}{\mathbb R}
\newcommand{\eps}{\varepsilon}
\newcommand{\To}{\longrightarrow}
\newcommand{\BX}{\mathbf{B}(X)}
\newcommand{\A}{\mathcal{A}}

\title{\bf  ADM submanifolds, SL normal bundles examples}
\author{\bf Doan The Hieu\\
\sl  University of Hue, 32 Le Loi, Hue, Vietnam\\
\sl deltic@dng.vnn.vn; dthehieu@yahoo.com} \maketitle


\begin{abstract}
It is showed  that many examples of AMD submanifolds of higher dimensions come from SL normal bundles. A symmetry property of SL submanifolds and Bj\"{o}rling type problem for SL normal bundles are also mentioned.
\end{abstract}
\vskip 1cm
\section{Introduction}

Throughout this paper, AMD means {\sl ``area-minimizing under diffeomorphisms leaving the boundary fixed''} and SL means  {\sl``special Lagrangian''} for short.

A submanifold $M\subset \Bbb R^n$   is said to be
AMD if
 $$Vol(M)\le
Vol(\varphi(M)),$$ for any diffeomorphism $\varphi$ of $\Bbb
R^n$ leaving the boundary $\partial M$ fixed.

 Some first examples of non flat sheets meeting along multiple curvers at an equal angle, that are AMD, were given in \cite{hi}. The non flat sheets are all calibrated surfaces of dimension two in $\Bbb R^4.$ The key idea is: try to find a suitable calibration $w,$ a corresponding calibrated submanifold $S$ and a plane $P$ of codimension 2 such that

1.  The intersection of $S$ and a hyperplane containing $P$ is a subset of $P$.

2. The resultings $w_i=R_{(i\alpha,P)}(w),\ i=0,1,2,\ldots, k-1$ have vanishing sum $\sum_{i=0}^{k-1}w_i=0$; where $R_{(i\alpha, P)}$ is the rotation about $P$ of  angle $i\alpha,\ \ k\alpha=2\pi.$
 
When all objects are found, we can choose a piecie of $S,$ denote by $S_0,$ that have  a part of the boundary in $P.$ The union of $S_i= R_{(i\alpha,P)}(S_0),$ \ $\cup_{i=0}^kS_i,$ that have a singular edge in $P$ can be proved to be AMD by using Stoke's theorem (see the example in Section 4) or Theorem ~\ref{t:main}.

In this paper, we show that one can choose the calibration $w$, the plane $P$ and the calibrated surface $S$  by a SL calibration, a complex plane of real codimention two and a suitable SL normal bundle. Also, a symmetry property of SL submanifolds and a Bj\"{o}rling type problem for SL normal bundles are mentioned.


\section{AMD sufficient condition}
In order to be easier for the readers, we recall some definitions and facts about AMD surmanifolds (see \cite {hi}).
\par 
Let $\{C_i\}_{i\in I}$ be a set of calibrated submanifolds of
dimension $m$ in $\Bbb R^n  (m<n)$ 
and $\{w_i\}_{i\in I}$ be the
set of correspondent calibrations. That means for each $i\in I,
w_i$ calibrates $C_i$ with a suitable orientation.  
Note that if $\omega_i$ calibrates $C_i$, then $-\omega_i$ calibrates $C_i$ with opposite orientation. Depending on a chosen orientation on $C_i$ we have the corespondent calibration to be $\omega_i$ or $-\omega_i.$

Let  $\Sigma\subset\Bbb R^n$ be a set satisfies the following conditions:

(i)\ $\Sigma\subset \cup_{i\in I}C_i,$

(ii)\ the set $E =\Sigma\cap(C_i\cap C_j)$ is of dimension $m-1$ for every $i,j\in I, i\not =j.$

We call each $F_i=\Sigma\cap C_i$ a face, each $E$ a singular
edge, the union of all singular edges $ E$ the singular set $S$,
the closure of $\partial F_i\sim S$ the boundary edge of $\Sigma$
in $F_i,$ the union $\cup_{i\in I}( \partial F_i\sim S)$ the boundary $\partial \Sigma$ of $\Sigma.$

Suppose $\{E_j\}_{j\in J}$ is the set of all singular edges and $\{F_i\}_{i\in I}$ is the set of all faces of $\Sigma.$  Denote
 $$I_{E_j}= \{i : F_i\supset E_j \}\subset I,$$
 $$J_{F_i}=\{j :  E_j\subset F_i \}\subset J.$$

\begin{thm}[Theorem 2.1 in \cite{hi}]\label{t:main}
Let $\Sigma$ be a set defined as above. Suppose that every singular edge $E_j$ lies on the boundary $\partial F_i, \forall i\in I_{E_j}$ and for each $E_j$ we can choose suitable orientations on $F_i, \forall i\in I_{E_j},$ such that:

(i) the  orientations on $F_i, \forall i\in I_{E_j}$ determine the same orientation on $E_j,$

(ii) the corespondent calibrations have vanishing sum.

 \noindent Then $\Sigma$ is
area-minimizing under diffeomorphisms leaving $\partial\Sigma$
fixed.
\end{thm}

\begin{cor}[Corollary 2.2 in \cite{hi}]\label{t:main2} 
Let $\Sigma$ be a polyhedral set.  Then $\Sigma$ is
area-minimizing under diffeomorphisms leaving $\partial\Sigma$
fixed if and only if  $\Sigma$   satisfies the assumptions in the Theorem    ~\ref{t:main}.
\end{cor}
\section{Normal bundle of an austere manifold}

It is well known that, for a $m$-submanifold $M$ of $\Bbb R^n$ its normal bundles $\nu(M) $ in $\Bbb C^n=\Bbb R^n\oplus\Bbb R^n$ is  a Lagrangian submanifold of $\Bbb R^n\oplus\Bbb R^n\cong \Bbb C^n$ with respect to the sympletic structure $w=dx_1\wedge dy_1+dx_2\wedge dy_2+\ldots+dx_n\wedge dy_n.$ Lawson and Harvey studied the case when  $\nu(M)$ become special Lagragian and proved the following. 
\begin{thm}[ ~\cite{hala}, Theorem 3.11, p. 102,]\label{hala}
The normal bundle $\nu(M)$ is special Lagrangian with phase $i^{n-m}$ if and only if $M$ is austere.
\end{thm}
A submanifold is called austere if its principal curvartures associated to any normal field is invariant under multiplication by -1.  Austere submanifolds were studied  by Bryant (see \cite{bry})  and Dajczer and Frorit (see \cite {daflo}). A minimal surfaces (of dimension two) is, of course, austere.
\par
Borisenko generalized this result by adding a harmonic function to  the case when $M$ is a minimal surface in $\Bbb R^3.$ 
\par
Let $M$ be a regular minimal surface in $\Bbb R^3$ and $\rho$ be a harmonic function on $M$. Let $p\in M$ and $(x,y)$ be local coordinates in a neighborhood of $p.$ Let
$$n(p)=\frac{r_x\times r_y}{|r_x\times r_y|}$$
be the normal vector and
$$\tau(p)=\frac{\rho_xr_y - \rho_yr_x}{|r_x\times r_y|}.$$

\begin{thm}[ ~\cite{bo}, Theorem 1]\label{bo}
The following submanifold in $\Bbb C^3=\Bbb R^3\oplus \Bbb R^3$
$$N=\{(p, \tau(p)\times n(p)+tn(p)\in \Bbb C^3: \ p\in M, t\in \Bbb R\}$$
is special Lagrangian. 
\end{thm}

\section{SL normal bundles examples}
In this section, we focus on constructing some AMD examples   from SL normal bundles. For details about SL normal bundles we refer the readers to \cite{bo}  and \cite{hala}. 
\par
For some AMD examples about polyhedral sets and the first ones about non flat sheets meeting along multiple curvers with an equal angle, see \cite{hi}. In order to illustrate, we begin by a simple example with a short proof. The proof is an application of Stoke's theorem  as that of the Fundamental theorem of calibrations (see \cite{hala}). 

Consider the following complex curve in the complex plane $\Bbb C^2\equiv \Bbb R^4$ with the
standard complex structure $J, Je_1=e_3; Je_2=e_4:$
\begin{eqnarray}
S&=&\{(z,w)\in \Bbb C^2 : z=w^2\}\nonumber\\
&=&\{(x_1, x_2, x_3, x_4) \in \Bbb R^4: \ \ x_2=x_1^2-x_3^2;\ \ x_4=2x_1x_3\}.\nonumber
\end{eqnarray}
\begin{eqnarray}
S_0&=S\cap \{(x_1, x_2, x_3, x_4) \in \Bbb R^4: \ \ x_1^2+x_2^2+x_3^2+x_4^2\le 1; x_3\ge 0\}.\nonumber
\end{eqnarray}  
By the Fundamental Theorem of calibrations, $S_0$ is calibrated by the Kahler form  $\omega_0$     and therefore is area-minimizing.

For $i=0,1,2,\ldots, k-1,$ let $R_i$ be the rotations of angles $\frac {2i\pi} k$ 
about the plane $\langle e_1, e_2\rangle.$  Obviously, $R_i\omega_0:=\omega_i, $ are Kahler forms corresponding to complex structures $R_iJ$ and calibrate $R_iS_0:=S_i$  for all $ i=0,1,2,\ldots, k-1,$ respectively. Direct computation shows that 
$$w_0+w_1+w_2+\ldots+w_{k-1}=0.$$
These complex curves $S_0, S_1, S_2,\ldots, S_{k-1}$  can be viewed as  two-dimensional,
area-minimizing (real) surfaces in $\Bbb R^4.$  Their boundaries contain
the common curve
$$K=\{x_1=x_2^2, \ x_3=x_4=0, x_1^2+x_2^2\le 1\},$$
that lies on the plane $\langle e_1, e_2\rangle.$

Let $\Sigma =\bigcup^{k-1}_{i=0} S_i.$ Then, Stoke's theorem yields:
  $$ \begin{aligned}
Vol(\Sigma)&=\Sigma_{i=0}^{k-1} Vol(S_i)=\Sigma_{i=0}^{k-1} \int_{S_i}w_i\\
&= \int_{S_1-S_0}w_1+\int_{S_2-S_0}w_2+\ldots+\int_{S_{k-1}-S_0}w_{k-1}\\
&= \int_{\varphi(S_1)-\varphi(S_0)}w_1+\int_{\varphi(S_2)-\varphi(S_0)}w_2+\ldots+\int_{\varphi(S_{k-1})-\varphi(S_0)}w_{k-1}\\
&= \Sigma_{i=0}^{k-1} \int_{\varphi(S_i)}w_i\le\Sigma_{i=0}^{k-1} Vol(\varphi(S_i))=Vol(\varphi(\Sigma)),\\
\end{aligned} 
$$
where  $\varphi$ is a diffeomorphism leaving $\partial \Sigma=\cup\partial S_i\sim K$ fixed.
 
As in the statement in the introduction and the above example, to construct AMD submanifolds from calibrated submanifolds, first we should find a suitable calibration and a plane of codimension two and second, a suitable calibrated submanifold.   
The following lemma shows that SL calibrations and complex planes of real codimension two satify the first requirement. 

Let $w=\text{Re}(dz_1\wedge\ldots\wedge dz_n)$ be the SL calibration on $\Bbb C^n$ with standard complex structure,  $P$ be a complex plane of real codimension two. Denote by $R_{(\alpha,P)}$ the rotation about the plane $P$ of angle $\alpha,\ k\alpha=2\pi,\ k\in \Bbb N\backslash\{0\}, \ 0<\alpha\le\pi.$ We have
\begin{lem}\label{spe}
Forms $w_i=R_{(i\alpha,P)}(w), i=0,1,2,\ldots k-1$ are all SL and have the vanishing sum
$$\sum_{i=0}^kw_k=0.$$
\end{lem}

{\bf Proof.}\ Obviously, each $R_{(i\alpha,P)}, i=0,1,2,\ldots k-1,$  belongs to $SU(n),$ and hence $w_i=R_{(i\alpha,P)}(w), i=0,1,2,\ldots k-1$ are all SL.

Now assume that $\{e_1, e_2\}$ is the real orthonormal basic of $P^\bot,$ the orthogonal comlement of $P.$ Without lost of generality, we can assume $w$ is of the following
$$w=e_1\wedge\varphi +e_2\wedge\psi,$$
where $\varphi=\text{Re}(dz_2\wedge\ldots\wedge dz_n)$ and $\psi=\text{Im} (dz_2\wedge\ldots\wedge dz_n).$ Direct computation shows that
$$w_i=(\cos i\alpha e_1+\sin i\alpha e_2)\wedge\varphi+(-\sin i\alpha e_1+\cos i\alpha e_2)\wedge\psi,$$ 
and hence
$$\sum_{i=0}^kw_k=0.$$
As the next step, we will show that there exits many SL submanifolds satisfying the second requirement.
The intersection of a such submanifold with a suitable part of hypersphere yields a  piece of SL submanifolds that have some parts of boundary lying in planes of real codimention two. Rotations of this about the planes (of suitable angles) give us examples of AMD manifolds by vitue of theorem ~\ref{t:main} or the proof of the above example. 

The following theorem shows that such SL submanifolds have a ``symmetry property." 

Let $S$ be a SL n-submanifold in $\Bbb C^n\equiv \Bbb R^{2n}$ that have a part of boundary $L\subset \partial S$ in a complex plane $P$ of real codimension two. Denote $S^*$ be the image of $S$ under the reflection $f$ along $P$ (this can be seen as the rotation about $P$ of angle $\pi$).

\begin{thm}[Symmetry of SL submanifolds]\label{sym}
$S\cup -S^*$ is a SL submanifold and hence area-minimizing.
\end{thm}
{\bf Proof.} Suppose $S$ is calibrated by SL-calibration $\varphi$.  Since $P$ is  a complex plane, $f$ belongs to $SU(n).$ It is easy to see that, $f(\varphi)= -\varphi$ and  $S^*$ is calibrated by $f(\varphi)=-\varphi.$  Thus,  $S\cup -S^*$ is calibrated by $\varphi$ and hence is a SL submanifold.

The following theorem shows that many SL 3-submanifolds that have a $4$-planar submanifold can be constructed in $\Bbb C^3.$
\begin{thm} [A Bj\"{o}rling type problem for SL normal bundles]
Let $\gamma_1,\gamma_2:I\longrightarrow \Bbb R^3$ be analytic curvers in $\Bbb R^3.$ Suppose that $\gamma_1'(t).\gamma_2(t)=0, \gamma_2(t)\not=0, \forall t\in I.$ Then  there exists a SL normal bundle in $\Bbb C^3\equiv \Bbb R^3\oplus\Bbb R^3$ contains the complex curve $\gamma_1(t)+i\gamma_2(t).$
\end{thm}
{\bf Proof.} First, there exists the solution $S$ to Bj\"{o}rling problem  for the curve $\gamma_1$ and normal field $\frac{\gamma_2(t)}{\|\gamma_2(t)\|}.$ This is a minimal surface that contains $\gamma_1$ and accepts $\frac{\gamma_2(t)}{\|\gamma_2(t)\|}$ as normal field. Since a minimal surface is austere,  the normal bundle $\nu(S)$ is a SL submanifold by vitue of Theorem ~\ref{hala}. Clearly, $\nu(S)$ contains the complex curve $\gamma_1+i\gamma_2.$ 

\par
The last part of this section is to present some explicit AMD examples with computations.
\par
{\bf Example 1.} Let $M$ be a minimal surface in $\Bbb R^3.$ If a plane meets $M$ orthogonally  everywhere along the intersection which is a curve, then that curve is called a  planar symmetry line of $M.$ If the curve is a straight line then it is called  straight symmetry line. Most of minimal surfaces have planar symmetry lines and straight symmetry lines. By the solution Bj\"{o}rling problem applying for a planar curve, we can construct  many such  minimal surfaces.
\par
Now, let $M$ be a minimal surface with a  planar symmetry lines $\gamma$ and $P$ be the plane containing $\gamma$ and meeting the surface $M$orthogonally. It is easy to see that, all normal vector of $M$ along $\gamma$ are in $P.$ Consider the intersection 
$$M_0=\nu(M)\cap S^+$$
where $S^+$ is a closed half-sphere, whose planar boundary lies in the hyperplane
$$\{(x, y)\in\Bbb R^3\oplus\Bbb R^3 :\ x\in P, y\in \Bbb R^3\}.$$
Since M is austere, $M_0$ is SL by Theorem ~\ref{hala}. We can see $M_0$ has a part of boundary that lies in 4-plane 
$$P+iP=\{(x,y) \in\Bbb R^3\oplus\Bbb R^3 : \ x\in P, y=\nu(x)\in \Bbb R^3\}.$$
This is the surface (of real dimention two)
$$\nu(\gamma)=\{(x,y) \in\Bbb R^3\oplus\Bbb R^3 :  |x|^2+|y|^2\le 1; \ x\in \gamma\cap S^+, y=\nu(x)\in \Bbb R^3\}.$$
Of course, $P+iP$ is complex 4-plane in $\Bbb C^3,$ and hence rotations about $P+iP$ of suitable angles give us 3-dimensional examples of AMD submanifolds by Lemma ~\ref{spe} and ~\ref{t:main}. 
\par
Below is an explicit example.
For $M$ let us take the catenoid
$$M=\{(u,\cosh u\cos v, \cosh u\sin v)\in \Bbb R^3\}.$$
Direct computations yields
$$n(u,v)=(\cosh u\sinh u, -\cosh u \cos v, -\cosh u\sin v).$$
and
$$\nu(M)=\{(u,\cosh u\cos v, \cosh u\sin v,w\cosh u\sinh u, -w\cosh u \cos v, -w\cosh u\sin v)\}$$
The intersection
$$M_0=\nu(M)\cap\{(x_1, x_2,x_3,y_1,y_2,y_3)\in \Bbb R^3\oplus \Bbb R^3: \sum x_i^2+\sum y_i^2\le 1; x_3\ge 0\}$$
has a 4-dimensional plannar part of boundary; and
$$M_1=\nu(M)\cap\{(x_1, x_2,0,y_1,y_2,y_3)\in \Bbb R^3\oplus \Bbb R^3:\sum x_i^2+\sum y_i^2\le 1; x_3\ge 0; x_2\ge 0\}$$
has two 4-dimensional plannar parts of boundary.

Rotations about planes $\{x_3=0; y_3=0\}$ and  $\{x_2=0; y_2=0\}$ of suitable angles, yield examples of 3-dimentional AMD submanifolds in $\Bbb C^3$ with one singular edge or many singular edges.

{\bf Example 2.}
As in example 1, but minimal surface $M$ replaces by an austere $m$-submanifold with a hyperplane meeting orthogonally along the intersection which is a $(m-1)$-submanifold. Below is an interesting example with infinite hyperplanes meeting the submanifold orthogonally.

 Consider the Clifford torus
$$T=\{(\cos u, \sin u, \cos v, \sin v) : u,v\in \Bbb R\}\subset\Bbb R^4.$$
It is well known that, $T$ is minimal in $S^3$ and hence austere. The cone over $T$
$$M={\cal C}(T)=\{(w\cos u, w\sin u, w\cos v, w\sin v) : u,v, w\in \Bbb R\}$$
is easily proved autere in $\Bbb R^4.$ 

A computation shows that
$$N=(-w^2\cos u\sin^2v,-w^2\sin u\sin^2v,-w^2\cos v\cos^2u,-w^2\sin v\cos^2v).$$ And we can verify that following hyperplanes
$$ax_1+bx_2=0;\  a,b\in \Bbb R;\ a^2+b^2\not=0$$
and 
$$cx_3+dx_4=0;\  c,d\in \Bbb R;\ c^2+d^2\not=0$$
meeting $T$ orthogonally.

{\bf Example 3.} We can use Theorem ~\ref{bo} to construct another examples. 
Below is an explicit example.
 $M$ is the catenoid
$$M=\{(u,\cosh u\cos v, \cosh u\sin v)\in \Bbb R^3\}$$
and let $\rho$ be  a harmonic function, for example
$$\rho=\sinh u\cos v.$$
Direct computations yields
$$n(u,v)=\frac 1{\cosh u}(\sinh u, -\cos v, -\sin v),$$
and
$$\tau (u,v)\times n(u,v)=\frac 1{\cosh u}(\cos v,\sin u, 0).$$
Consider
$$\nu(M)=\{((u,\cosh u\cos v, \cosh u\sin v,\frac{\cos v}{\cosh u}+w\tanh u, \tanh u- w\frac{\cos v}{\cosh u},-w\frac{\sin v}{\cosh u})\}$$
The intersection
$$M_0=\nu(M)\cap\{(x_1, x_2,x_3,y_1,y_2,y_3)\in \Bbb R^3\oplus \Bbb R^3: \sum x_i^2+\sum y_i^2\le 1; x_3\ge 0\}$$
has a 4-dimensional plannar part of boundary.

Rotations about planes $\{x_3=0; y_3=0\}$  of suitable angles, yields another examples of 3-dimentional AMD submanifolds in $\Bbb C^3.$ 
\par

\end{document}